# A General Class of Estimators of Population Median Using Two Auxiliary Variables in Double Sampling

**Jack Allen[1], Housila P. Singh[2], Sarjinder Singh[3], Florentin Smarandache[4]**

[1] School of Accounting and Finance, Griffith University, Australia
[2] School of Studies in Statistics, Vikram University, Ujjain - 456 010 (M. P.), India
[3] Department of Mathematics and Statistics, University of Saskatchewan, Canada
[4] Department of Mathematics, University of New Mexico, Gallup, USA

**Abstract:**
In this paper we have suggested two classes of estimators for population median $M_Y$ of the study character Y using information on two auxiliary characters X and Z in double sampling. It has been shown that the suggested classes of estimators are more efficient than the one suggested by Singh *et al* (2001). Estimators based on estimated optimum values have been also considered with their properties. The optimum values of the first phase and second phase sample sizes are also obtained for the fixed cost of survey.

1. INTRODUCTION

In survey sampling, statisticians often come across the study of variables which have highly skewed distributions, such as income, expenditure etc. In such situations, the estimation of median deserves special attention. Kuk and Mak (1989) are the first to introduce the estimation of population median of the study variate Y using auxiliary information in survey sampling. Francisco and Fuller (1991) have also considered the problem of estimation of the median as part of the estimation of a finite population distribution function. Later Singh *et al* (2001) have dealt extensively with the problem of estimation of median using auxiliary information on an auxiliary variate in two phase sampling.

Consider a finite population U={1,2,…,i,...,N}. Let Y and X be the variable for study and auxiliary variable, taking values $Y_i$ and $X_i$ respectively for the i-th unit. When the two variables are strongly related but no information is available on the population median $M_X$ of X, we seek to estimate the population median $M_Y$ of Y from a sample $S_m$, obtained through a two-phase selection. Permitting simple random sampling without replacement (SRSWOR) design in each phase, the two-phase sampling scheme will be as follows:

(i)  The first phase sample $S_n(S_n \subset U)$ of fixed size n is drawn to observe only X in order to furnish an estimate of $M_X$.

(ii) Given $S_n$, the second phase sample $S_m(S_m \subset S_n)$ of fixed size m is drawn to observe Y only.

Assuming that the median $M_X$ of the variable X is known, Kuk and Mak (1989) suggested a ratio estimator for the population median $M_Y$ of Y as

$$\hat{M}_1 = \hat{M}_Y \frac{M_X}{\hat{M}_X} \tag{1.1}$$

where $\hat{M}_Y$ and $\hat{M}_X$ are the sample estimators of $M_Y$ and $M_X$ respectively based on a sample $S_m$ of size m. Suppose that $y_{(1)}, y_{(2)}, \ldots, y_{(m)}$ are the y values of sample units in ascending order. Further, let t be an integer such that $Y_{(t)} \leq M_Y \leq Y_{(t+1)}$ and let p=t/m be the proportion of Y, values in the sample that are less than or equal to the median value $M_Y$, an unknown population parameter. If $\hat{p}$ is a predictor of p, the sample median $\hat{M}_Y$ can be written in terms of quantities as $\hat{Q}_Y(\hat{p})$ where $\hat{p} = 0.5$. Kuk and Mak (1989) define a matrix of proportions ($P_{ij}(x,y)$) as

|  | Y ≤ $M_Y$ | Y > $M_Y$ | Total |
|---|---|---|---|
| X ≤ $M_X$ | $P_{11}(x,y)$ | $P_{21}(x,y)$ | $P_{\cdot 1}(x,y)$ |
| X > $M_X$ | $P_{12}(x,y)$ | $P_{22}(x,y)$ | $P_{\cdot 2}(x,y)$ |
| Total | $P_{1\cdot}(x,y)$ | $P_{2\cdot}(x,y)$ | 1 |

and a position estimator of $M_y$ given by

$$\hat{M}_Y^{(p)} = \hat{Q}_Y(\hat{p}_Y) \qquad (1.2)$$

where 
$$\hat{p}_Y = \frac{1}{m}\left(\frac{m_x \hat{p}_{11}(x,y)}{\hat{p}_{.1}(x,y)} + \frac{(m-m_x)\hat{p}_{12}(x,y)}{\hat{p}_{.2}(x,y)}\right)$$
$$\approx 2\left(\frac{m_x \hat{p}_{11}(x,y) + (m-m_x)\hat{p}_{12}(x,y)}{m}\right)$$

with $\hat{p}_{ij}(x,y)$ being the sample analogues of the $P_{ij}(x,y)$ obtained from the population and $m_x$ the number of units in $S_m$ with $X \leq M_X$.

Let $\tilde{F}_{YA}(y)$ and $\tilde{F}_{YB}(y)$ denote the proportion of units in the sample $S_m$ with $X \leq M_X$, and $X > M_X$, respectively that have Y values less than or equal to y. Then for estimating $M_Y$, Kuk and Mak (1989) suggested the 'stratification estimator' as

$$\hat{M}_Y^{(St)} = \inf\{y : \tilde{F}_Y^{(y)} \geq 0.5\} \qquad (1.3)$$

where $\hat{F}_Y(y) \cong \frac{1}{2}\left[\tilde{F}_{YA}^{(y)} + \tilde{F}_{YB}^{(y)}\right]$

It is to be noted that the estimators defined in (1.1), (1.2) and (1.3) are based on prior knowledge of the median $M_X$ of the auxiliary character X. In many situations of practical importance the population median $M_X$ of X may not be known. This led Singh *et al* (2001) to discuss the problem of estimating the population median $M_Y$ in double sampling and suggested an analogous ratio estimator as

$$\hat{M}_{1d} = \hat{M}_Y \frac{\hat{M}_X^1}{\hat{M}_X} \qquad (1.4)$$

where $\hat{M}_X^1$ is sample median based on first phase sample $S_n$.

Sometimes even if $M_X$ is unknown, information on a second auxiliary variable Z, closely related to X but compared X remotely related to Y, is available on all units of the population. This type of situation has been briefly discussed by, among others, Chand (1975), Kiregyera (1980, 84), Srivenkataramana and Tracy (1989), Sahoo and Sahoo (1993) and Singh (1993). Let $M_Z$ be the known population median of Z. Defining

$$e_0 = \left(\frac{\hat{M}_Y}{M_Y} - 1\right), e_1 = \left(\frac{\hat{M}_X}{M_X} - 1\right), e_2 = \left(\frac{\hat{M}_X^1}{M_X} - 1\right), e_3\left(\frac{\hat{M}_Z}{M_Z} - 1\right) \text{ and } e_4 = \left(\frac{\hat{M}_Z^1}{M_Z} - 1\right)$$

such that $E(e_k) \cong 0$ and $|e_k| < 1$ for $k=0,1,2,3$; where $\hat{M}_2$ and $\hat{M}_2^1$ are the sample median estimators based on second phase sample $S_m$ and first phase sample $S_n$. Let us define the following two new matrices as

|  | $Z \leq M_Z$ | $Z > M_Z$ | Total |
|---|---|---|---|
| $X \leq M_X$ | $P_{11}(x,z)$ | $P_{21}(x,z)$ | $P_{\cdot 1}(x,z)$ |
| $X > M_X$ | $P_{12}(x,z)$ | $P_{22}(x,z)$ | $P_{\cdot 2}(x,z)$ |
| Total | $P_{1\cdot}(x,z)$ | $P_{2\cdot}(x,z)$ | 1 |

and

|  | $Z \leq M_Z$ | $Z > M_Z$ | Total |
|---|---|---|---|
| $Y \leq M_Y$ | $P_{11}(y,z)$ | $P_{21}(y,z)$ | $P_{\cdot 1}(y,z)$ |
| $Y > M_Y$ | $P_{12}(y,z)$ | $P_{22}(y,z)$ | $P_{\cdot 2}(y,z)$ |
| Total | $P_{1\cdot}(y,z)$ | $P_{2\cdot}(y,z)$ | 1 |

Using results given in the Appendix-1, to the first order of approximation, we have

$$E(e_0^2) = \left(\frac{N-m}{N}\right)(4m)^{-1}\{M_Y f_Y(M_Y)\}^{-2},$$

$$E(e_1^2) = \left(\frac{N-m}{N}\right)(4m)^{-1}\{M_X f_X(M_X)\}^{-2},$$

$$E(e_2^2) = \left(\frac{N-n}{N}\right)(4n)^{-1}\{M_X f_X(M_X)\}^{-2},$$

$$E(e_3^2) = \left(\frac{N-m}{N}\right)(4m)^{-1}\{M_Z f_Z(M_Z)\}^{-2},$$

$$E(e_4^2) = \left(\frac{N-n}{N}\right)(4n)^{-1}\{M_Z f_Z(M_Z)\}^{-2},$$

$$E(e_0 e_1) = \left(\frac{N-m}{N}\right)(4m)^{-1}\{4P_{11}(x,y)-1\}\{M_X M_Y f_X(M_X) f_Y(M_Y)\}^{-1},$$

$$E(e_0 e_2) = \left(\frac{N-n}{N}\right)(4n)^{-1}\{4P_{11}(x,y)-1\}\{M_X M_Y f_X(M_X) f_Y(M_Y)\}^{-1},$$

$$E(e_0 e_3) = \left(\frac{N-m}{N}\right)(4m)^{-1}\{4P_{11}(y,z)-1\}\{M_Y M_Z f_Y(M_Y) f_Z(M_Z)\}^{-1},$$

$$E(e_0 e_4) = \left(\frac{N-n}{N}\right)(4n)^{-1}\{4P_{11}(y,z)-1\}\{M_Y M_Z f_Y(M_Y) f_Z(M_Z)\}^{-1},$$

$$E(e_1 e_2) = \left(\frac{N-n}{N}\right)(4n)^{-1}\{M_X f_X(M_X)\}^{-2},$$

$$E(e_1 e_3) = \left(\frac{N-m}{N}\right)(4m)^{-1}\{4P_{11}(x,z)-1\}\{M_X M_Z f_X(M_X) f_Z(M_Z)\}^{-1},$$

$$E(e_1 e_4) = \left(\frac{N-n}{N}\right)(4n)^{-1}\{4P_{11}(x,z)-1\}\{M_X M_Z f_X(M_X) f_Z(M_Z)\}^{-1},$$

$$E(e_2 e_3) = \left(\frac{N-n}{N}\right)(4n)^{-1}\{4P_{11}(x,z)-1\}\{M_X M_Z f_X(M_X) f_Z(M_Z)\}^{-1},$$

$$E(e_2e_4) = \left(\frac{N-n}{N}\right)(4n)^{-1}\{4P_{11}(x,z)-1\}\{M_XM_Zf_X(M_X)f_Z(M_Z)\}^{-1},$$

$$E(e_3e_4) = \left(\frac{N-n}{N}\right)(4n)^{-1}(f_Z(M_Z)M_Z)^{-2}$$

where it is assumed that as $N \to \infty$ the distribution of the trivariate variable $(X,Y,Z)$ approaches a continuous distribution with marginal densities $f_X(x)$, $f_Y(y)$ and $f_Z(z)$ for X, Y and Z respectively. This assumption holds in particular under a superpopulation model framework, treating the values of $(X, Y, Z)$ in the population as a realization of N independent observations from a continuous distribution. We also assume that $f_Y(M_Y)$, $f_X(M_X)$ and $f_Z(M_Z)$ are positive.

Under these conditions, the sample median $\hat{M}_Y$ is consistent and asymptotically normal (Gross, 1980) with mean $M_Y$ and variance

$$\left(\frac{N-m}{N}\right)(4m)^{-1}\{f_Y(M_Y)\}^{-2}$$

In this paper we have suggested a class of estimators for $M_Y$ using information on two auxiliary variables X and Z in double sampling and analyzes its properties.

## 2. SUGGESTED CLASS OF ESTIMATORS

Motivated by Srivastava (1971), we suggest a class of estimators of $M_Y$ of Y as

$$g = \{\hat{M}_Y^{(g)} : \hat{M}_Y^{(g)} = M_Y g(u,v)\} \qquad (2.1)$$

where $u = \dfrac{\hat{M}_X}{\hat{M}_X^1}, v = \dfrac{\hat{M}_Z^1}{\hat{M}_Z}$ and g(u,v) is a function of u and v such that g(1,1)=1 and such that it satisfies the following conditions.

1. Whatever be the samples ($S_n$ and $S_m$) chosen, let (u,v) assume values in a closed convex sub-space, P, of the two dimensional real space containing the point (1,1).

2. The function g(u,v) is continuous in P, such that g(1,1)=1.

3. The first and second order partial derivatives of g(u,v) exist and are also continuous in P.

Expanding g(u,v) about the point (1,1) in a second order Taylor's series and taking expectations, it is found that

$$E(\hat{M}_Y^{(g)}) = M_Y + 0(n^{-1})$$

so the bias is of order $n^{-1}$.

Using a first order Taylor's series expansion around the point (1,1) and noting that g(1,1)=1, we have

$$\hat{M}_Y^{(g)} \cong M_Y[1 + e_0 + (e_1 - e_2)g_1(1,1) + e_4 g_2(1,1) + 0(n^{-1})]$$

or

$$(M_Y^{(g)} - M_Y) \cong M_Y[e_0 + (e_1 - e_2)g_1(1,1) + e_4 g_2(1,1)] \quad (2.2)$$

where $g_1(1,1)$ and $g_2(1,1)$ denote first order partial derivatives of g(u,v) with respect to u and v respectively around the point (1,1).

Squaring both sides in (2.2) and then taking expectations, we get the variance of $\hat{M}_Y^{(g)}$ to the first degree of approximation, as

$$Var(\hat{M}_Y^{(g)}) = \frac{1}{4(f_Y(M_Y))^2}\left[\left(\frac{1}{m} - \frac{1}{N}\right) + \left(\frac{1}{m} - \frac{1}{n}\right)A + \left(\frac{1}{n} - \frac{1}{N}\right)B\right], \quad (2.3)$$

where

$$A = \left(\frac{M_Y f_Y(M_Y)}{M_X f_X(M_X)}\right)g_1(1,1)\left[\left(\frac{M_Y f_Y(M_Y)}{M_X f_X(M_X)}\right)g_1(1,1) + 2(4P_{11}(x,y) - 1)\right] \quad (2.4)$$

$$B = \left(\frac{M_Y f_Y(M_Y)}{M_Z f_Z(M_Z)}\right)g_Z(1,1)\left[\left(\frac{M_Y f_Y(M_Y)}{M_Z f_Z(M_Z)}\right)g_2(1,1) + 2(4P_{11}(y,z) - 1)\right] \quad (2.5)$$

The variance of $\hat{M}_Y^{(g)}$ in (2.3) is minimized for

$$g_1(1,1) = -\left(\frac{M_X f_X(M_X)}{M_Y f_Y(M_Y)}\right)(4P_{11}(x,y) - 1)$$

$$g_2(1,1) = -\left(\frac{M_Z f_Z(M_Z)}{M_Y f_Y(M_Y)}\right)(4P_{11}(y,z) - 1)$$

(2.6)

Thus the resulting (minimum) variance of $M_Y^{(g)}$ is given by

$$\text{min. Var}(\hat{M}_Y^{(g)}) = \frac{1}{4(f_Y(M_Y))^2}\left[\left(\frac{1}{m} - \frac{1}{N}\right) - \left(\frac{1}{m} - \frac{1}{n}\right)(4P_{11}(x,y) - 1)^2 - \left(\frac{1}{n} - \frac{1}{N}\right)(4P_{11}(y,z) - 1)\right]$$

(2.7)

Now, we proved the following theorem.

Theorem 2.1 - Up to terms of order $n^{-1}$,

$$\text{Var}(\hat{M}_Y^g) \geq \frac{1}{4(f_y(M_Y))^2}\left[\left(\frac{1}{m}-\frac{1}{N}\right)-\left(\frac{1}{m}-\frac{1}{n}\right)(4P_{11}(x,y)-1)^2-\left(\frac{1}{n}-\frac{1}{N}\right)(4P_{11}(y,z)-1)^2\right]$$

with equality holding if

$$g_1(1,1) = -\left(\frac{M_x f_x(M_x)}{M_Y f_Y(M_Y)}\right)(4P_{11}(x,y)-1)$$

$$g_2(1,1) = -\left(\frac{M_z f_z(M_z)}{M_Y f_Y(M_Y)}\right)(4P_{11}(y,z)-1)$$

It is interesting to note that the lower bound of the variance of $\hat{M}_y^{(g)}$ at (2.1) is the variance of the linear regression estimator

$$\hat{M}_Y^{(l)} = \hat{M}_Y + \hat{d}_1(\hat{M}_X^1 - \hat{M}_X) + \hat{d}_2(M_Z - \hat{M}_Z^1) \quad (2.8)$$

where

$$\hat{d}_1 = \frac{\hat{f}_X(\hat{M}_x)}{\hat{f}_Y(\hat{M}_y)}(4\hat{p}_{11}(x,y)-1),$$

$$\hat{d}_2 = \frac{\hat{f}_Z(\hat{M}_Z)}{\hat{f}_Y(\hat{M}_Y)}(4\hat{p}_{11}(y,z)-1),$$

with $\hat{p}_{11}(x,y)$ and $\hat{p}_{11}(y,z)$ being the sample analogues of the $p_{11}(x,y)$ and $p_{11}(y,z)$ respectively and $\hat{f}_Y(\hat{M}_Y), \hat{f}_X(M_X)$ and $\hat{f}_Z(M_Z)$ can be obtained by following Silverman (1986).

Any parametric function g(u,v) satisfying the conditions (1), (2) and (3) can generate an asymptotically acceptable estimator. The class of such estimators are large. The following simple functions g(u,v) give even estimators of the class

$$g^{(1)}(u,v) = u^\alpha v^\beta, \quad g^{(2)}(u,v) = \frac{1+\alpha(u-1)}{1-\beta(v-1)},$$

$$g^{(3)}(u,v) = 1+\alpha(u-1)+\beta(v-1), \quad g^{(4)}(u,v) = \{1-\alpha(u-1)-\beta(v-1)\}^{-1}$$

$$g^{(5)}(u,v) = w_1 u^\alpha + w_2 v^\beta, \quad w_1 + w_2 = 1$$

$$g^{(6)}(u,v) = \alpha u + (1-\alpha)v^\beta, \quad g^{(7)}(u,v) = \exp\{\alpha(u-1)+\beta(v-1)\}$$

Let the seven estimators generated by $g^{(i)}(u,v)$ be denoted by $\hat{M}_{Yi}^{(g)} = \hat{M}_Y g^{(i)}(u,v), (i=1 \text{ to } 7)$. It is easily seen that the optimum values of the parameters $\alpha, \beta, w_i(i-1,2)$ are given by the right hand sides of (2.6).

## 3. A WIDER CLASS OF ESTIMATORS

The class of estimators (2.1) does not include the estimator

$$\hat{M}_{Yd} = \hat{M}_Y + d_1(\hat{M}_X^1 - M_X) + d_2(M_Z - \hat{M}_Z^1), (d_1, d_2)$$

being constants.

However, it is easily shown that if we consider a class of estimators wider than (2.1), defined by

$$\hat{M}_Y^{(G)} = G_1(\hat{M}_Y, u, v) \tag{3.1}$$

of $M_Y$, where $G(\cdot)$ is a function of $\hat{M}_Y$, u and v such that $G(M_Y,1,1) = M_Y$ and $G_1(M_Y,1,1) = 1$. $G_1(M_Y,1,1)$ denoting the first partial derivative of $G(\cdot)$ with respect to $\hat{M}_Y$.

Proceeding as in Section 2 it is easily seen that the bias of $\hat{M}_Y^{(G)}$ is of the order $n^{-1}$ and up to this order of terms, the variance of $\hat{M}_Y^{(G)}$ is given by

$$\mathrm{Var}(\hat{M}_Y^{(G)}) = \frac{1}{4(f_Y(M_Y))^2} \left[ \left(\frac{1}{m} - \frac{1}{N}\right) + \left(\frac{1}{m} - \frac{1}{n}\right)\left(\frac{f_Y(M_Y)}{M_X f_X(M_X)}\right) \right.$$

$$G_2(M_Y,1,1) \left\{ \left(\frac{f_Y(M_Y)}{M_X f_X(M_X)}\right) G_2(M_Y,1,1) + 2(4P_{11}(x,y) - 1) \right\}$$

$$\left. + \left(\frac{1}{n} - \frac{1}{N}\right) \frac{f_Y(M_Y)}{f_Z(M_Z) M_Z} \left\{ \left(\frac{f_Y(M_Y)}{M_Z f_Z(M_Z)}\right) G_3(M_Y,1,1) + 2(4P_{11}(y,z) - 1) \right\} \right]$$

(3.2)

where $G_2(M_Y1,1)$ and $G_3(M_Y1,1)$ denote the first partial derivatives of u and v respectively around the point $(M_Y,(1,1))$.

The variance of $\hat{M}_Y^{(G)}$ is minimized for

$$G_2(M_Y,1,1) = -\left(\frac{M_X f_X(M_X)}{f_Y(M_Y)}\right)(4P_{11}(x,y) - 1)$$

$$G_3(M_Y,1,1) = -\left(\frac{M_Z f_Z(M_Z)}{f_Y(M_Y)}\right)(4P_{11}(y,z) - 1)$$

(3.3)

Substitution of (3.3) in (3.2) yields the minimum variance of $\hat{M}_Y^{(G)}$ as

$$\text{min. Var}\!\left(\hat{M}_Y^{(G)}\right) = \frac{1}{4(f_Y(M_Y))^2}\left[\left(\frac{1}{m}-\frac{1}{N}\right)-\left(\frac{1}{m}-\frac{1}{n}\right)(4P_{11}(x,y)-1)^2 -\left(\frac{1}{n}-\frac{1}{N}\right)(4P_{11}(y,z)-1)^2\right]$$
$$= \text{min.Var}\!\left(\hat{M}_Y^{(g)}\right)$$

(3.4)

Thus we established the following theorem. Theorem 3.1 - Up to terms of order $n^{-1}$,

$$\text{Var}\!\left(\hat{M}_Y^{(G)}\right) \geq \frac{1}{4(f_Y(M_Y))^2}\left[\left(\frac{1}{m}-\frac{1}{N}\right)-\left(\frac{1}{m}-\frac{1}{n}\right)(4P_{11}(x,y)-1)^2 -\left(\frac{1}{n}-\frac{1}{N}\right)(4P_{11}(y,z)-1)^2\right]$$

with equality holding if

$$G_2(M_Y,1,1) = -\left(\frac{f_x(M_X)M_X}{f_Y(M_Y)}\right)(4P_{11}(x,y)-1)$$

$$G_3(M_Y,1,1) = -\left(\frac{M_Z f_Z(M_Z)}{f_Y(M_Y)}\right)(4P_{11}(y,z)-1)$$

If the information on second auxiliary variable z is not used, then the class of estimators $\hat{M}_Y^{(G)}$ reduces to the class of estimators of $M_Y$ as

$$\hat{M}_Y^{(H)} = H\!\left(\hat{M}_Y, u\right) \tag{3.5}$$

where $H\!\left(\hat{M}_Y,u\right)$ is a function of $\left(\hat{M}_Y,u\right)$ such that $H(M_Y,1) = M_Y$ and $H_1(M_Y,1)=1$,

$H_1(M_Y,1) = \dfrac{\partial H(\cdot)}{\partial \hat{M}_Y}\bigg|_{(M_Y,1)}$. The estimator $\hat{M}_Y^{(H)}$ is reported by Singh *et al* (2001).

The minimum variance of $\hat{M}_Y^{(H)}$ to the first degree of approximation is given by

$$\text{min.Var}\!\left(\hat{M}_Y^{(H)}\right) = \frac{1}{4(f_Y(M_Y))^2}\left[\left(\frac{1}{m}-\frac{1}{N}\right)-\left(\frac{1}{m}-\frac{1}{n}\right)(4P_{11}(x,y)-1)^2\right] \tag{3.6}$$

From (3.4) and (3.6) we have

$$\text{minVar}\!\left(\hat{M}_Y^{(H)}\right) - \text{min.Var}\!\left(\hat{M}_Y^{(G)}\right) = \left(\frac{1}{n}-\frac{1}{N}\right)\frac{1}{4(f_Y(M_Y))^2}(4P_{11}(y,z)-1)^2 \tag{3.7}$$

which is always positive. Thus the proposed class of estimators $\hat{M}_Y^{(G)}$ is more efficient than the estimator $\hat{M}_Y^{(H)}$ considered by Singh *et al* (2001).

## 4. ESTIMATOR BASED ON ESTIMATED OPTIMUM VALUES

We denote

$$\alpha_1 = \frac{M_X f_X(M_X)}{M_Y f_Y(M_Y)}(4P_{11}(x,y)-1)$$
$$\alpha_2 = \frac{M_Z f_Z(M_Z)}{M_Y f_Y(M_Y)}(4P_{11}(y,z)-1)$$

(4.1)

In practice the optimum values of $g_1(1,1)(=-\alpha_1)$ and $g_2(1,1)(=-\alpha_2)$ are not known. Then we use to find out their sample estimates from the data at hand. Estimators of optimum value of $g_1(1,1)$ and $g_2(1,1)$ are given as

$$\hat{g}_1(1,1) = -\hat{\alpha}_1$$
$$\hat{g}_2(1,1) = -\hat{\alpha}_2$$

(4.2)

where

$$\hat{\alpha}_1 = \frac{\hat{M}_X \hat{f}_X(\hat{M}_X)}{\hat{M}_Y \hat{f}_Y(\hat{M}_Y)}(4\hat{p}_{11}(x,y)-1)$$
$$\hat{\alpha}_2 = \frac{\hat{M}_Z \hat{f}_Z(\hat{M}_Z)}{\hat{M}_Y \hat{f}_Y(\hat{M}_Y)}(4p_{11}(y,z)-1)$$

(4.3)

Now following the procedure discussed in Singh and Singh (19xx) and Srivastava and Jhajj (1983), we define the following class of estimators of $M_Y$ (based on estimated optimum) as

$$\hat{M}_Y^{(g*)} = \hat{M}_Y g*(u,v,\hat{\alpha}_1,\hat{\alpha}_2)$$

(4.4)

where $g*(\cdot)$ is a function of $(u,v,\hat{\alpha}_1,\hat{\alpha}_2)$ such that

$g*(1,1,\alpha_1\alpha_2) = 1$

$g_1^*(1,1,\alpha_1,\alpha_2) = \left.\frac{\partial g*(\cdot)}{\partial u}\right|_{(1,1,\alpha_1,\alpha_2)} = -\alpha_1$

$g_2^*(1,1,\alpha_1,\alpha_2) = \left.\frac{\partial g*(\cdot)}{\partial v}\right|_{(1,1,\alpha_1,\alpha_2)} = -\alpha_2$

$g_3^*(1,1,\alpha_1,\alpha_2) = \left.\frac{\partial g*(\cdot)}{\partial \hat{\alpha}_1}\right|_{(1,1,\alpha_1,\alpha_2)} = 0$

$g_4^*(1,1,\alpha_1,\alpha_2) = \left.\frac{\partial g*(\cdot)}{\partial \hat{\alpha}_2}\right|_{(1,1,\alpha_1,\alpha_2)} = 0$

and such that it satisfies the following conditions:

1. Whatever be the samples ($S_n$ and $S_m$) chosen, let $u, v, \hat{\alpha}_1, \hat{\alpha}_2$ assume values in a closed convex sub-space, S, of the four dimensional real space containing the point $(1,1,\alpha_1,\alpha_2)$.

2. The function $g^*(u, v, \alpha_1, \alpha_2)$ continuous in S.

3. The first and second order partial derivatives of $g^*(u, v, \hat{\alpha}_1, \hat{\alpha}_2)$ exst. and are also continuous in S.

Under the above conditions, it can be shown that

$$E(\hat{M}_Y^{(g^*)}) = M_Y + 0(n^{-1})$$

and to the first degree of approximation, the variance of $\hat{M}_Y^{(g^*)}$ is given by

$$\text{Var}(\hat{M}_Y^{(g^*)}) = \min.\text{Var}(\hat{M}_Y^g) \qquad (4.5)$$

where $\min.\text{Var}(\hat{M}_Y^{(g)})$ is given in (2.7).

A wider class of estimators of $M_Y$ based on estimated optimum values is defined by

$$\hat{M}_Y^{(G^*)} = G^*(\hat{M}_Y, u, v, \hat{\alpha}_1^*, \hat{\alpha}_2^*) \qquad (4.6)$$

where

$$\hat{\alpha}_1^* = \frac{\hat{M}_X \hat{f}_X(\hat{M}_X)}{\hat{f}_Y(\hat{M}_Y)}(4\hat{p}_{11}(x, y) - 1)$$

$$\hat{\alpha}_2^* = \frac{\hat{M}_Z \hat{f}_Z(\hat{M}_Z)}{\hat{f}_Y(\hat{M}_Y)}(4\hat{p}_{11}(y, z) - 1) \qquad (4.7)$$

are the estimates of

$$\alpha_1^* = \frac{M_X f_X(M_X)}{f_Y(M_Y)}(4P_{11}(x, y) - 1)$$

$$\alpha_2^* = \frac{M_Z f_Z(M_Z)}{f_Y(M_Y)}(4P_{11}(y, z) - 1) \qquad (4.8)$$

and $G^*(\cdot)$ is a function of $(\hat{M}_Y, u, v, \alpha_1^*, \hat{\alpha}_2^*)$ such that

$$G^*(M_Y,1,1,\alpha_1^*,\alpha_2^*) = M_Y$$

$$G_1^*(M_Y,1,1,\alpha_1^*,\alpha_2^*) = \left.\frac{\partial G^*(\cdot)}{\partial \hat{M}_Y}\right|_{(M_Y,1,1,\alpha_1^*,\alpha_2^*)} = 1$$

$$G_2^*(M_Y,1,1,\alpha_1^*,\alpha_2^*) = \left.\frac{\partial G^*(\cdot)}{\partial u}\right|_{(M_Y,1,1,\alpha_1^*,\alpha_2^*)} = -\alpha_1^*$$

$$G_3^*(M_Y,1,1,\alpha_1^*,\alpha_2^*) = \left.\frac{\partial G^*(\cdot)}{\partial v}\right|_{(M_Y,1,1,\alpha_1^*,\alpha_2^*)} = -\alpha_2^*$$

$$G_4^*(M_Y,1,1,\alpha_1^*,\alpha_2^*) = \left.\frac{\partial G^*(\cdot)}{\partial \hat{\alpha}_1^*}\right|_{(M_Y,1,1,\alpha_1^*,\alpha_2^*)} = 0$$

$$G_5^*(M_Y,1,1,\alpha_1^*,\alpha_2^*) = \left.\frac{\partial G^*(\cdot)}{\partial \hat{\alpha}_2^*}\right|_{(M_Y,1,1,\alpha_1^*,\alpha_2^*)} = 0$$

Under these conditions it can be easily shown that

$$E(\hat{M}_Y^{(G^*)}) = M_Y + 0(n^{-1})$$

and to the first degree of approximation, the variance of $\hat{M}_Y^{(G^*)}$ is given by

$$\mathrm{Var}(\hat{M}_Y^{G^*}) = \min.\mathrm{Var}(\hat{M}_Y^{(G)}) \qquad (4.9)$$

where $\min.\mathrm{Var}(\hat{M}_Y^{G})$ is given in (3.4).

It is to be mentioned that a large number of estimators can be generated from the classes $\hat{M}_Y^{(g^*)}$ and $\hat{M}_Y^{(G^*)}$ based on estimated optimum values.

## 5. EFFICIENCY OF THE SUGGESTED CLASS OF ESTIMATORS FOR FIXED COST

The appropriate estimator based on on single-phase sampling without using any auxiliary variable is $\hat{M}_Y$, whose variance is given by

$$\mathrm{Var}(\hat{M}_Y) = \left(\frac{1}{m} - \frac{1}{N}\right)\frac{1}{4(f_Y(M_Y))^2} \qquad (5.1)$$

In case when we do not use any auxiliary character then the cost function is of the form $C_0 - mC_1$, where $C_0$ and $C_1$ are total cost and cost per unit of collecting information on the character Y.

The optimum value of the variance for the fixed cost $C_0$ is given by

$$\text{Opt.}\left[\text{Var}(\hat{M}_Y)\right] = V_0\left(\frac{G}{C_0} - \frac{1}{N}\right) \tag{5.2}$$

where

$$V_0 \frac{1}{4(f_Y(M_Y))^2} \tag{5.3}$$

When we use one auxiliary character X then the cost function is given by

$$C_0 = Gm + C_2 n, \tag{5.4}$$

where $C_2$ is the cost per unit of collecting information on the auxiliary character Z.

The optimum sample sizes under (5.4) for which the minimum variance of $\hat{M}_Y^{(H)}$ is optimum, are

$$m_{opt} = \frac{C_0 \sqrt{(V_0 - V_1)/C_1}}{\left[\sqrt{(V_0 - V_1)C_1} + \sqrt{V_1 C_2}\right]} \tag{5.5}$$

$$n_{opt} = \frac{C_0 \sqrt{V_1/C_2}}{\left[\sqrt{(V_0 - V_1)C_1} + \sqrt{V_1 C_2}\right]}$$

where $V_1 = V_0(4P_{11}(x,y) - 1)^2$.

Putting these optimum values of m and n in the minimum variance expression of $\hat{M}_Y^{(H)}$ in (3.6), we get the optimum $\min.\text{Var}(\hat{M}_Y^{(H)})$ as

$$\text{Opt.}\left[\min.\text{Var}(\hat{M}_Y^{(H)})\right] = \left[\frac{\left(\sqrt{(V_0 - V_1)C_1} + \sqrt{V_1 C_2}\right)^2}{C_0} - \frac{V_0}{N}\right] \tag{5.7}$$

Similarly, when we use an additional character Z then the cost function is given by

$$C_0 = C_1 m + (C_2 + C_3)n \tag{5.8}$$

where $C_3$ is the cost per unit of collecting information on character Z.

It is assumed that $C_1 > C_2 > C_3$. The optimum values of m and n for fixed cost $C_0$ which minimizes the minimum variance of $\hat{M}_Y^{(g)}\left(or \hat{M}_Y^{(G)}\right)$ (2.7) (or (3.4)) are given by

$$m_{opt} = \frac{C_0 \sqrt{(V_0 - V_1)/C_1}}{\left[\sqrt{(V_0 - V_1)C_1} + \sqrt{(C_2 + C_3)(V_1 - V_2)}\right]} \quad (5.9)$$

$$n_{opt} = \frac{C_0 \sqrt{(V_1 - V_2)/C_2 + C_3}}{\left[\sqrt{(V_0 - V_1)C_1} + \sqrt{(C_2 + C_3)(V_1 - V_2)}\right]} \quad (5.10)$$

where $V_2 = V_0(4P_{11}(y,z)-1)^2$.

The optimum variance of $\hat{M}_Y^{(g)}\left(or \hat{M}_Y^{(G)}\right)$ corresponding to optimal two-phase sampling strategy is

$$\mathrm{Opt}\left[\min.\mathrm{Var}\left(\hat{M}_Y^{(g)}\right) or \min.\mathrm{Var}\left(\hat{M}_Y^{(G)}\right)\right] = \left[\frac{[\sqrt{(V_0 - V_1)C_1} + \sqrt{(C_2 + C_3)(V_1 - V_2)}]^2}{C_0} - \frac{V_2}{N}\right] \quad (5.11)$$

Assuming large N, the proposed two phase sampling strategy would be profitable over single phase sampling so long as

$$\left[\mathrm{Opt}.\mathrm{Var}\left(\hat{M}_Y\right)\right] > \mathrm{Opt}.\left[\min.\mathrm{Var}\left(\hat{M}_Y^{(g)}\right) or \min.\mathrm{Var}\left(\hat{M}_Y^{(G)}\right)\right]$$

i.e. $\dfrac{C_2 + C_3}{C_1} < \left[\dfrac{\sqrt{V_0} - \sqrt{V_0 - V_1}}{\sqrt{V_1 - V_2}}\right] \quad (5.12)$

When N is large, the proposed two phase sampling is more efficient than that Singh *et al* (2001) strategy if

$$\mathrm{Opt}\left[\min.\mathrm{Var}\left(\hat{M}_Y^{(g)}\right) or \min.\mathrm{Var}\left(\hat{M}_Y^{(G)}\right)\right] < \mathrm{Opt}\left[\min.\mathrm{Var}\left(\hat{M}_Y^{(H)}\right)\right]$$

i.e. $\dfrac{C_2 + C_3}{C_1} < \dfrac{V_1}{V_1 - V_2} \quad (5.13)$

6. GENERALIZED CLASS OF ESTIMATORS

We suggest a class of estimators of $M_Y$ as

$$\Im = \{\hat{M}_Y^{(F)} : \hat{M}_Y^{(F)} = F(\hat{M}_Y, u, v, w)\} \qquad (6.1)$$

where $u = \hat{M}_X / \hat{M}'_X$, $v = \hat{M}'_Z / M_Z$, $w = \hat{M}_Z / M_Z$ and the function $F(\cdot)$ assumes a value in a bounded closed convex subset $W \subset \Re_4$, which contains the point $(M_Y,1,1,1)=T$ and is such that $F(T)=M_Y \Rightarrow F_1(T)=1$, $F_1(T)$ denoting the first order partial derivative of $F(\cdot)$ with respect to $\hat{M}_Y$ around the point $T=(M_Y,1,1,1)$. Using a first order Taylor's series expansion around the point T, we get

$$\hat{M}_Y^{(F)} = F(T) + (\hat{M}_Y = M_Y) F_1(T) + (u-1)F_2(T) + (v-1)F_3(T) + (w-1)F_4(T) + 0(n^{-1}) \qquad (6.2)$$

where $F_2(T)$, $F_3(T)$ and $F_4(T)$ denote the first order partial derivatives of $F(\hat{M}_Y, u, v, w)$ with respect to $u$, $v$ and $w$ around the point T respectively. Under the assumption that $F(T)=M_Y$ and $F_1(T)=1$, we have the following theorem.

Theorem 6.1. Any estimator in $\Im$ is asymptotically unbiased and normal.

Proof: Following Kuk and Mak (1989), let $P_Y$, $P_X$ and $P_Z$ denote the proportion of Y, X and Z values respectively for which $Y \leq M_Y$, $X \leq M_X$ and $Z \leq M_Z$; then we have

$$\hat{M}_Y - M_Y = \frac{1}{2 f_Y(M_Y)}(1 - 2P_Y) + 0_p\left(n^{-\frac{1}{2}}\right)$$

$$\hat{M}_X - M_X = \frac{1}{2 f_X(M_X)}(1 - 2P_X) + 0_p\left(n^{-\frac{1}{2}}\right)$$

$$\hat{M}'_X - M_X = \frac{1}{2 f_X(M_X)}(1 - 2P_X) + 0_p\left(n^{-\frac{1}{2}}\right)$$

$$\hat{M}_Z - M_Z = \frac{1}{2 f_Z(M_Z)}(1 - 2P_Z) + 0_p\left(n^{-\frac{1}{2}}\right)$$

and

$$\hat{M}'_Z - M_Z = \frac{1}{2 f_Z(M_Z)}(1 - 2P_Z) + 0_p\left(n^{-\frac{1}{2}}\right)$$

Using these expressions in (6.2), we get the required results.

Expression (6.2) can be rewritten as

$$\hat{M}_Y^{(F)} - M_Y \cong (\hat{M}_Y - M_Y) + (u-1)F_2(T) + (v-1)F_3(T) + (w-1)F_4(T)$$

or

$$\hat{M}_Y^{(F)} - M_Y \cong M_Y e_0 + (e_1 - e_2) F_2(T) + e_4 F_3(T) + e_3 F_4(T) \tag{6.3}$$

Squaring both sides of (6.3) and then taking expectation, we get the variance of $\hat{M}_Y^{(F)}$ to the first degree of approximation, as

$$\operatorname{Var}(\hat{M}_Y^{(F)}) = \frac{1}{4(f_Y(M_Y))^2}\left[\left(\frac{1}{m} - \frac{1}{N}\right)A_1 + \left(\frac{1}{m} - \frac{1}{n}\right)A_2 + \left(\frac{1}{n} - \frac{1}{N}\right)A_3\right], \tag{6.4}$$

where

$$A_1 = \left[1 + \left(\frac{f_Y(M_Y)}{M_Z f_Z(M_Z)}\right)^2 F_4^{\,2}(T) + 2(4P_{11}(y,z)-1)\left(\frac{f_Y(M_Y)}{M_Z f_Z(M_Z)}\right)F_4(T)\right]$$

$$A_2 = \left(\frac{f_Y(M_Y)}{M_X f_X(M_X)}\right)\left[\left(\frac{f_Y(M_Y)}{M_X f_X(M_X)}\right)F_2^{\,2}(T) + 2(4P_{11}(x,y)-1)F_2(T) \right.$$
$$\left. + 2(4P_{11}(x,z)-1)\left(\frac{f_Y(M_Y)}{M_Z f_Z(M_Z)}\right)F_2(T)F_4(T)\right]$$

$$A_3 = \left(\frac{f_Y(M_Y)}{M_Z f_Z(M_Z)}\right)\left[\left(\frac{f_Y(M_Y)}{M_Z f_Z(M_Z)}\right)F_3^{\,2}(T) + 2(4P_{11}(y,z)-1)F_3(T) \right.$$
$$\left. + 2\left(\frac{f_Y(M_Y)}{M_Z f_Z(M_Z)}\right)F_3(T)F_4(T)\right]$$

The $\operatorname{Var}(\hat{M}_Y^{(F)})$ at (6.4) is minimized for

$$F_2(T) = -\frac{[(4P_{11}(x,y)-1)-(4P_{11}(x,z)-1)(4P_{11}(y,z)-1)]}{[1-(4P_{11}(x,z)-1)^2]} \cdot \frac{M_X f_X(M_X)}{f_Y(M_Y)}$$
$$= -a_2 \text{ (say)}$$

$$\tag{6.5}$$

$$F_3(T) = -\frac{(4P_{11}(x,z)-1)[(4P_{11}(x,y)-1)-(4P_{11}(y,z)-1)(4P_{11}(x,z)-1)]}{[1-(4P_{11}(x,z)-1)^2]} \cdot \frac{M_Z f_Z(M_Z)}{f_Y(M_Y)}$$
$$= -a_2 \text{ (say)}$$

$$F_4(T) = -\frac{[(4P_{11}(y,z)-1)-(4P_{11}(x,y)-1)(4P_{11}(x,z)-1)]}{[1-(4P_{11}(x,z)-1)^2]} \cdot \frac{M_z f_z(M_z)}{f_Y(M_Y)}.$$

$$= -a_3 \text{ (say)}$$

Thus the resulting (minimum) variance of $\hat{M}_Y^{(F)}$ is given by

$$\min \text{Var}(\hat{M}_Y^{(F)}) = \frac{1}{4(f_Y(M_Y))^2} \left[ \begin{array}{c} \left(\frac{1}{m}-\frac{1}{N}\right) - \left(\frac{1}{m}-\frac{1}{n}\right)\left\{\frac{D^2}{1-(4P_{11}(x,z)-1)^2} + (4P_{11}(x,y)-1)\right\}^2 \\ -\left(\frac{1}{n}-\frac{1}{N}\right)(4P_{11}(y,z)-1)^2 \end{array} \right]$$

$$= \min.\text{Var}(\hat{M}_Y^{(G)}) - \left(\frac{1}{m}-\frac{1}{n}\right)\frac{1}{4(f_Y(M_Y))^2} \frac{D^2}{1-[4P_{11}(x,z)-1^2]}$$

(6.6)

where

$$D = [(4P_{11}(y,z)-1)-(4P_{11}(x,y)-1)(4P_{11}(x.z)-1)] \tag{6.7}$$

and $\min.\text{Var}(\hat{M}_Y^{(G)})$ is given in (3.4)

Expression (6.6) clearly indicates that the proposed class of estimators $\hat{M}_Y^{(F)}$ is more efficient than the class of estimator $\hat{M}_Y^{(G)}$ or $(\hat{M}_Y^{(g)})$ and hence the class of estimators $\hat{M}_Y^{(H)}$ suggested by Singh *et al* (2001) and the estimator $\hat{M}_Y$ at its optimum conditions.

The estimator based on estimated optimum values is defined by

$$p^* = \{\hat{M}_Y^{(F*)} : \hat{M}_Y^{F*} = F^*(\hat{M}_Y, u, v, w, \hat{a}_1, \hat{a}_2, \hat{a}_3)\} \tag{6.8}$$

where

$$\hat{a}_1 = \frac{[(4\hat{p}_{11}(x,y)-1)-(4\hat{p}_{11}(x,z)-1)(4\hat{p}_{11}(y,z)-1)]}{[1-(4\hat{p}_{11}(x,z)-1)^2]} \frac{\hat{M}_x \hat{f}_x(\hat{M}_x)}{\hat{f}_Y(\hat{M}_Y)}.$$

$$\hat{a}_2 = \frac{(4\hat{p}_{11}(x,z)-1)[(4\hat{p}_{11}(x,y)-1)-(4\hat{p}_{11}(y,z)-1)(4\hat{p}_{11}(x,z)-1)]}{[1-(4\hat{p}_{11}(x,z)-1)^2]} \frac{\hat{M}_z \hat{f}_z(\hat{M}_z)}{\hat{f}_Y(\hat{M}_Y)}.$$

$$a_3 = \frac{[(4\hat{p}_{11}(y,z)-1)-(4\hat{p}_{11}(x,y)-1)(4\hat{p}_{11}(x,z)-1)]}{[1-(4\hat{p}_{11}(x,z)-1)^2]} \frac{\hat{M}_Z \hat{f}_Z(\hat{M}_Z)}{\hat{f}_Y(\hat{M}_Y)}.$$

(6.9)

are the sample estimates of $a_1$, $a_2$ and $a_3$ given in (6.5) respectively, F*(·) is a function of $(\hat{M}_Y, u, v, w, \hat{a}_1, \hat{a}_2, \hat{a}_3)$ such that

$F^*(T^*) = M_Y$

$$\Rightarrow F_1^*(T^*) = \left.\frac{\partial F^*(\cdot)}{\partial \hat{M}_Y}\right|_{T^*} = 1$$

$$F_2^*(T^*) = \left.\frac{\partial F^*(\cdot)}{\partial u}\right|_{T^*} = -a_1$$

$$F_3^*(T^*) = \left.\frac{\partial F^*(\cdot)}{\partial v}\right|_{T^*} = -a_2$$

$$F_4^*(T^*) = \left.\frac{\partial F^*(\cdot)}{\partial w}\right|_{T^*} = -a_3$$

$$F_5^*(T^*) = \left.\frac{\partial F^*(\cdot)}{\partial \hat{a}_1}\right|_{T^*} = 0$$

$$F_6^*(T^*) = \left.\frac{\partial F^*(\cdot)}{\partial \hat{a}_2}\right|_{T^*} = 0$$

$$F_7^*(T^*) = \left.\frac{\partial F^*(\cdot)}{\partial \hat{a}_3}\right|_{T^*} = 0$$

where $T^* = (M_Y, 1, 1, 1, a_1, a_2, a_3)$

Under these conditions it can easily be shown that

$$E(\hat{M}_Y^{(F^*)}) = M_Y + 0(n^{-1})$$

and to the first degree of approximation, the variance of $\hat{M}_Y^{(F^*)}$ is given by

$$\text{Var}(\hat{M}_Y^{(F^*)}) = \min.\text{Var}(\hat{M}_Y^F)$$

(6.10)

where $\min.\text{Var}\left(\hat{M}_Y^{(F)}\right)$ is given in (6.6).

Under the cost function (5.8), the optimum values of m and n which minimizes the minimum variance of $\hat{M}_Y^{(F)}$ is (6.6) are given by

$$m_{opt} = \frac{C_0\sqrt{(V_0 - V_1 - V_3)/C_1}}{[\sqrt{(V_0 - V_1 - V_3)C_1} + \sqrt{(V_1 - V_2 - V_3)(C_2 + C_3)}]} \tag{6.11}$$

$$n_{opt} = \frac{C_0\sqrt{(V_1 - V_2 - V_3)/C_2}}{[\sqrt{(V_0 - V_1 - V_3)C_1} + \sqrt{(V_1 - V_2 + V_3)(C_2 + C_3)}]}$$

where

$$V_3 = \frac{D^2 V_0}{\left[1 - (4P_{11}(x,z) - 1)^2\right]} \tag{6.12}$$

for large N, the optimum value of $\min.\text{Var}\left(\hat{M}_Y^{(F)}\right)$ is given by

$$\text{Opt.}\left[\min.\text{Var}\left(\hat{M}_Y^{(F)}\right)\right] = \frac{\left[\sqrt{(V_0 - V_1 - V_3)C_1} + \sqrt{(V_1 - V_2 + V_3)(C_2 + C_3)}\right]}{C_0} \tag{6.13}$$

The proposed two-phase sampling strategy would be profitable over single phase-sampling so long as $\text{Opt.}\left[\text{Var}(\hat{M}_Y)\right] > \text{Opt.}\left[\min.\text{Var}(\hat{M}_Y^{(F)})\right]$

$$\text{i.e.} \quad \frac{C_2 + C_3}{C_1} < \left[\frac{\sqrt{V_0} - \sqrt{V_0 - V_1 - V_3}}{\sqrt{V_1 - V_2 + V_3}}\right]^2 \tag{6.14}$$

It follows from (5.7) and (6.13) that

$$\text{Opt.}\left[\min.\text{Var}\left(\hat{M}_Y^{(F)}\right)\right] < \text{Opt.}\left[\min.\text{Var}\left(\hat{M}_Y^H\right)\right]$$

$$\text{if} \quad \left(\frac{\sqrt{V_0 - V_1} - \sqrt{V_0 - V_1 - V_3}}{\sqrt{V_1 - V_2 + V_3}}\right) > \left[\sqrt{\frac{C_2 + C_3}{C_1}} - \sqrt{\frac{V_1}{(V_1 - V_2 + V_3)C_1}\frac{C_2}{C_1}}\right] \tag{6.15}$$

for large N.

Further we note from (5.11) and (6.13) that

$$\text{Opt.}\left[\min.\text{Var}\left(\hat{M}_Y^{(F)}\right)\right] < \text{Opt.}\left[\min.\text{Var}\left(\hat{M}_Y^{(g)} \text{or} \hat{M}_Y^G\right)\right]$$

$$\text{if} \quad \frac{C_2 + C_3}{C_1} < \left[ \frac{\sqrt{(V_0 - V_1)} - \sqrt{(V_0 - V_1 - V_3)}}{\sqrt{(V_1 - V_2 + V_3)} - \sqrt{V_1 - V_2}} \right]^2 \quad (6.16)$$